\documentclass{amsart}
\input{xy}
\xyoption{all}

\newcommand{\K}{\mathcal K}
\newcommand{\I}{\mathcal I}
\newcommand{\IR}{\mathbb R}
\newcommand{\Ra}{\Rightarrow}

\newcommand{\w}{\omega}
\newcommand{\dlim}{\varinjlim}

\newcommand{\Top}{\mathbf{Top}}
\newcommand{\Tych}{\mathbf{Tych}}

\newcommand{\Comp}{\mathbf{Comp}}

\newcommand{\C}{\mathbf{C}}
\newcommand{\ANE}{\mathsf{ANE}}
\newcommand{\IN}{\mathbb N}

\newcommand{\ANRsw}{\mathsf{ANR}(k_\omega)}

\newcommand{\Kw}{\mathbf K_\w}
\newcommand{\id}{\mathrm{id}}
\newcommand{\pr}{\mathrm{pr}}

\newtheorem{theorem}{Theorem}[section]

\newtheorem{lemma}[theorem]{Lemma}
\newtheorem{corollary}[theorem]{Corollary}

\theoremstyle{definition}
\newtheorem{definition}[theorem]{Definition}

\title[Free topological universal algebras]{Free topological universal algebras and absolute neighborhood retracts}
\author{Taras Banakh and Olena Hryniv}
\address{Ivan Franko national University of Lviv, Ukraine}
\email{t.o.banakh@gmail.com, olena\_hryniv@ukr.net}
\begin{document}
\begin{abstract} We prove that for a complete quasivariety $\K$ of topological $E$-algebras of countable discrete signature $E$ and each submetrizable  $\ANRsw$-space $X$ its free topological $E$-algebra $F_\K(X)$ in the class $\K$ is a submetrizable $\ANRsw$-space.
\end{abstract}

\maketitle

\section{Introduction}

In this paper we study the construction of a free topological universal algebra and show that this construction preserves the class of submetrizable $\ANRsw$-spaces.

To give a precise formulation of our main result, we need to recall some definitions related to topological universal algebras. For more detail information, see \cite{Chob1}, \cite{Chob2}, \cite{Chob3}.

\begin{definition} Let $(E_n)_{n\in\w}$ be a sequence of pairwise disjoint topological spaces. The topological sum $E=\bigoplus_{n\in\w}E_n$ is called {\em the continuous signature}. The signature is called {\em discrete} ({\em countable}) if so is the space $E$.

A  {\em  topological universal algebra of signature $E$} or briefly, a {\em topological $E$-algebra} is a topological space $X$ endowed with a family of continuous maps $e_{n,X}:E_n\times X^n\to X$, $n\in\w$.

A topological $E$-algebra $(X,\{e_{n,X}\}_{n\in\w})$ is called {\em Tychonov} if the underlying topological space $X$ is Tychonov.
\end{definition}

Homomorphisms between $E$-algebras are defined as follows.

\begin{definition}
A function $h:X\to Y$ between two topological $E$-algebras $(X,\{e_{n,X}\}_{n\in\w})$ and $(Y,\{e_{n,Y}\}_{n\in\w})$ is called a {\em $E$-homomorphism} if $$e_{n,Y}(z,h(x_1),\dots,h(x_n))=h(e_{n,X}(z,x_1,\dots,x_n))$$for any $n\in\w$, $z\in E_n$, and $x_1,\dots,x_n\in X$.

Such a function $h$ is called an {\em algebraic isomorphism} ({\em topological isomorphism}) if $h$ is bijective and both functions $h$ and $h^{-1}$ are (continuous) $E$-homomorphisms of the $E$-algebras.
\end{definition}

Next, we define some operations over $E$-algebras.

\begin{definition}
For topological $E$-algebras $X_\alpha$, $\alpha\in A$, the Tychonov product $X=\prod_{\alpha\in A}X_\alpha$ is a topological $E$-algebra endowed with the structure mappings $$e_{n,X}(z,x_1,\dots,x_n)=\big(e_{n,X_\alpha}(z,\pr_\alpha(x_1),\dots,\pr_\alpha(x_n))\big)_{\alpha\in A}$$
where $n\in\w$, $z\in E_n$, $x_1,\dots,x_n\in X$, and $\pr_\alpha:\prod_{\alpha\in A}X_\alpha\to X_\alpha$ is the $\alpha$-coordinate projection.
\end{definition}

\begin{definition} A subset $A\subset X$ of a topological $E$-algebra $(X,\{e_n\}_{n\in\w})$ is 
called a {\em subalgebra} if $e_n(E_n\times A^n)\subset A$ for all $n\in\w$.
\end{definition}

Since for any subalgebras $A_i\subset X$, $i\in\I$, of a topological $E$-algebra $X$ the intersection $A=\bigcap_{i\in\I}A_i$ is a subalgebra of $X$, for each subset $Z\subset X$ there is a minimal subalgebra $\langle Z\rangle$ of $X$ that contains $Z$. This is the subalgebra generated by the set $Z$. The structure of this subalgebra $\langle Z\rangle$ can be described as follows.

Given a subset $L\subset E$ and a subset $Z$ of a topological $E$-algebra $(X,\{e_n\}_{n\in\w})$, let $$
\begin{aligned}
&\langle Z\rangle^L_0=Z,\\
&\langle Z\rangle^L_{n+1}=\langle Z\rangle_{n}^L\cup \bigcup_{k\in\w}e_{k,X}\big((E_k\cap L)\times (\langle Z\rangle^L_{n})^k\big)\mbox{ for $n\in\w$, and }\\
&\langle Z\rangle_\w^L=\bigcup_{n\in\w}\langle Z\rangle^L_n.
\end{aligned}$$

By induction, one can check that for compact subspaces $L\subset E$ and $Z\subset X$ the subset $\langle Z\rangle^L_n$ of $X$ is compact for every $n\in\w$. Consequently, $\langle Z\rangle^L_\w$ is a $\sigma$-compact subset of $X$.

Writing the signature $E$ and the space $Z$ as the unions $E=\bigcup_{n\in\w}L_n$ and $Z=\bigcup_{n\in\w}Z_n$ of  non-decreasing sequences of subsets, we see that $$\langle Z\rangle=\bigcup_{n\in\w}\langle Z_n\rangle^{L_n}_n$$is the subalgebra of $X$, generated by $Z$. If the spaces $Z_n$ and $L_n$, $n\in\w$, are compact (finite), then each subset $\langle Z_n\rangle^{L_n}_n$, $n\in\w$, of $X$ is compact (finite) and hence the algebraic hull $\langle Z\rangle$ of $Z$ in $X$ is $\sigma$-compact (at most countable).

\begin{definition} A class $\K$ of topological $E$-algebras is called a {\em complete quasivariety} if
\begin{enumerate}
\item for each topological $E$-algebra $X\in\K$, each $E$-subalgebra of $X$ belongs to the class $\K$;
\item for any topological $E$-algebras $X_\alpha\in\K$, $\alpha\in A$, their  Tychonov product $\prod_{\alpha\in A}X_\alpha$ belongs to the class $\K$;
\item a Tychonov $E$-algebra belongs to $\K$ if it is algebraically isomorphic to a topological $E$-algebra $Y\in\K$.
\end{enumerate}
A complete quasivariety $\K$ is {\em non-trivial} if it contains a topological $E$-algebra $X$ that contains more that one point.
\end{definition} 

Finally, we recall the notion of a free topological $E$-algebra.

\begin{definition} Let $\K$ be a complete quasivariety of topological $E$-algebras. A {\em free topological $E$-algebra in $\K$} over a topological space $X$ is a pair $(F_\K(X),\eta)$ consisting of a topological $E$-algebra $F_\K(X)\in\K$ and a continuous map $\eta:X\to F_\K(X)$ such that for any continuous map $f:X\to Y$ to a topological $E$-algebra $Y\in\K$ there is a unique continuous $E$-homomorphism $h:F_\K(X)\to Y$ such that $f=h\circ \eta$.
\end{definition}

The construction $F_\K(X)$ of a free topological $E$-algebra has been intensively studied by M.M.Choban \cite{Chob2}, \cite{Chob3}. In particular, he proved that for each complete quasivariety $\K$ of topological $E$-algebras and any topological space $X$ a free topological $E$-algebra $(F_\K(X),\eta)$ exists and is unique up to a topological isomorphism. Also he proved the following important result, see \cite[2.4]{Chob2}:

\begin{theorem}[Choban]\label{choban1} If $\K$ is a non-trivial complete quasivariety of topological $E$-algebras, then for each Tychonov space $X$ the canonical map $\eta:X\to F_\K(X)$ is a topological embedding and $F_\K(X)$ coincides with the subalgebra $\langle\eta(X)\rangle$ generated by the image $\eta(X)$ of $X$ in $F(X,\K)$.
\end{theorem}

Since $\eta:X\to F_\K(X)$ is a topological embedding, we can identify a Tychonov space $X$ with its image $\eta(X)$ in $F_\K(X)$ and say that the free $E$-algebra $F_\K(X)$ is algebraically generated by $X$.

In fact, the construction of a free topological $E$-algebra $F_\K(X)$ determines a functor $F_\K:\Top\to \K$ from the category $\Top$ of topological spaces and their continuous maps to the category whose objects are topological $E$-algebras from the class $\K$ and morphisms are continuous $E$-homomorphisms.

In \cite{Chob1}--\cite{Chob3} a lot of attention  was paid to the problem of preservation of various topological properties by the functor $F_\K$. In particular, it was shown that the functor $F_\K$ preserves (submetrizable) $k_\w$-spaces provided the signature $E$ is a (submetrizable) $k_\w$-space, see \cite[4.1.2]{Chob3}.

A Hausdorff topological space $X$ is called {\em a $k_\w$-space} if $X=\dlim X_n$ is the {\em direct limit} of a non-decreasing sequence of compact subsets $(X_n)_{n\in\w}$ of $X$ in the sense that $X=\bigcup_{n\in\w}X_n$ and a subset $U\subset X$ is open if and only if $U\cap X_n$ is open in $X_n$ for each $n\in\w$. Such a sequence $(X_n)_{n\in\w}$ is called a {\em $k_\w$-sequence} for $X$.  

An {\em $s_\w$-space} is a direct limit $\dlim X_n$ of a $k_\w$-sequence $(X_n)_{n\in\w}$ consisting of second countable compact subspaces of $X$. It is easy to see that a $k_\w$-space $X$ is an $s_\w$-space if and only if it is {\em submetrizable} in the sense that $X$ admits a continuous metric.

\begin{theorem}[Choban]\label{choban2} Let $\K$ be a complete quasivariety of topological $E$-algebras  whose signature $E$ is a (submetrizable) $k_\w$-space. Then for each (submetrizable) $k_\w$-space $X$ the free topological $E$-algebra $F_\K(X)$ is a (submetrizable) $k_\w$-space. Moreover, if $E=\dlim L_n$ and $X=\dlim X_n$ for some $k_\w$-sequences $(L_n)_{n\in\w}$ and $(X_n)_{n\in\w}$,  then $(\langle \eta(X_n)\rangle^{L_n}_n)_{n\in\w}$ is a $k_\w$-sequence for $F_\K X$ and thus $F_\K X=\dlim\langle \eta(X_n)\rangle^{L_n}_n$. 
\end{theorem}

The principal result of this paper asserts that the functor $F_\K$ preserves $\ANRsw$-spaces.

\begin{definition} A $k_\w$-space $X$ is called an {\em absolute neighborhood retract in the class of $k_w$-spaces} (briefly, an $\ANRsw$) if $X$ is a neighborhood retract in each $k_\w$-space that contains $X$ as a closed subspace. 
\end{definition}

In Theorem~\ref{ANRkw} we shall show that a submetrizable $k_\w$-space $X$ is an $\ANRsw$-space if and only if $X$  each map $f:B\to X$ defined on a closed subspace of a (metrizable) compact space extends to a continuous map $\bar f:N(B)\to X$ defined on a neighborhood $N(B)$ of $B$ in $A$.

A topological space $X$ is called {\em compactly finite-dimensional} if each compact subset of $X$ is finite-dimensional. 

The following theorem is the main result of this paper.

\begin{theorem}\label{main}  If $\K$ is a complete quasivariety of topological $E$-algebras of countable discrete signature $E$, then for each submetrizable \textup{(}compactly finite-dimensional\textup{)} $\ANRsw$-space $X$ so is its free topological $E$-algebra $F_\K X$ in the quasivariety $\K$. 
\end{theorem}

In \cite{BH} this theorem will be applied for recognizing the topological structure of free topological inverse semigroups and groupoids.

\section{$\ANRsw$-spaces}

In this section we collect some information about $\ANRsw$-spaces. Such spaces are tightly connected with ANE-spaces.

Following \cite{Hu} we define a topological space $X$ to be an {\em absolute neighborhood extensor} for a class $\C$ of topological spaces (briefly, an $\ANE(\C)$-{\em space}) if each map $f:B\to X$ defined on a closed subspace $B$ of a topological space $C\in\C$ has a continuous extension $\bar f:N(B)\to X$ defined on some neighborhood $N(B)$ of $B$ in $C$. If always $f$ can be extended to the whole space $C$, then $X$ is called an {\em abslute extensor} for the class $\C$.

 By the Dugundji-Borsuk Theorem \cite{Dug}, \cite{Borsuk} each convex subset of a locally convex linear topological space, is an absolute extensor for the class of metrizable spaces. This theorem was generalized by Borges \cite{Borges} who proved that a convex subset of a locally convex space is an absolute extensor for the class of stratifiable spaces. This class contains all metrizable spaces and all submetrizable $k_\w$-spaces, and is closed with respect to many countable topological operations, see \cite{Borges}, \cite{Grue}. 

An important example of an $\ANRsw$-space is the space $$Q^\infty=\{(x_i)_{i\in\w}\in\IR^\infty:\sup_{i\in\w}|x_i|<\infty\}$$of bounded sequences, endowed with the direct limit topology $\dlim [-n,n]^\w$ generated by the $k_\w$-sequence $([-n,n]^\w)_{n\in\IN}$ consisting of the Hilbert cubes.
Being a locally convex linear topological space, $Q^\infty$ is an absolute extensor for the class of stratifiable spaces.

A topological space $X$ is called a {\em $Q^\infty$-manifold} if $X$ is Lindel\"of and each point $x\in X$ has a neighborhood homeomorphic to an open subset of $Q^\infty$. The theory of $Q^\infty$-manifolds was developed by K.Sakai \cite{Sak1}, \cite{Sak2} who established the following fundamental results:

\begin{theorem}[Characterization]\label{char} A topological space $X$ is homeomorphic to \textup{(}a manifold modeled on\textup{)} the space $Q^\infty$ if and only if $X$ is a submetrizable  $k_\w$-space such that each embedding $f:B\to X$ of a closed subset $B$ of a compact metrizable space $A$ can be extended to a topological embedding of \textup{(}an open neighborhood of $B$ in\textup{)} the space $A$ into $X$.
\end{theorem}

\begin{theorem}[Open Embedding]\label{opemb} Each $Q^\infty$-manifold is homeomorphic to an open subset of $Q^\infty$.
\end{theorem}

\begin{theorem}[Closed Embedding]\label{clemb} Each submetrizable $k_\w$-space is homeomorphic to a closed subspace of $Q^\infty$.
\end{theorem}

\begin{theorem}[Classification] Two $Q^\infty$-manifolds are homeomorphic if and only if they are homtopically equivalent.
\end{theorem}

\begin{theorem}[Triangulation]\label{triangle} Each $Q^\infty$-manifold $X$ is homeomorphic to $K\times Q^\infty$ for some countable locally finite simplicial complex $K$.
\end{theorem}

\begin{theorem}[ANR-Theorem]\label{ANR} For each submetrizable $\ANRsw$-space $X$ the product $X\times Q^\infty$ is a $Q^\infty$-manifold.
\end{theorem}

We shall use these theorems in the proof of the following (probably known as a folklor) characterization of submetrizable $\ANRsw$-spaces. 

\begin{theorem}\label{ANRkw} For a submetrizable $k_\w$-space $X$ the following conditions are equivalent:
\begin{enumerate}
\item $X$ is an $\ANRsw$-space;
\item $X$ is an ANE for the class of $k_\w$-spaces;
\item $X$ is an ANE for the class of compact metrizable spaces;
\item $X$ is an ANE for the class of stratifiable spaces;
\item $X$ is a retract of a $Q^\infty$-manifold.
\end{enumerate}
The equivalent conditions (1)-(5) hold if $X=\dlim X_n$ is the direct limit of a $k_\w$-sequence consisting of compact ANR's.
\end{theorem}

\begin{proof} $(1)\Ra(5)$ Assume that $X$ is an $\ANRsw$-space. By the Closed Embedding Theorem~\ref{clemb}, we can identify the submetrizable $k_\w$-space $X$ with a closed subspace of $Q^\infty$. Being an $\ANRsw$, $X$ is a retract of an open neighborhood $N(X)\subset Q^\infty$. Since $N(X)$ is a $Q^\infty$-manifold, $X$ is a retract of a $Q^\infty$-manifold.
\smallskip

$(5)\Ra(4)$ Assume that $X$ is a retract of a $Q^\infty$-manifold $M$. By the Open Embedding Theorem~\ref{opemb}, $M$ can be identified with an open subspace of $Q^\infty$. By the Borges' Theorem \cite{Borges}, the locally convex space $Q^\infty$ is an absolute extensor for the class of stratifiable spaces. Then the open subspace $M$ of $Q^\infty$ is an ANE for this class and so is its retract $X$.
\smallskip

The implication $(4)\Ra(3)$ is trivial since each metrizable space is stratifiable.
\smallskip

$(3)\Ra(2)$ Assume that $X$ is an ANE for the class of compact metrizable spaces. First we prove that $X$ is an ANE for the class of compact Hausdorff spaces. Let $f:B\to X$ be a continuous map defined on a closed subspace $B$ of a compact Hausdorff space $A$. Embed the compact space $A$ into a Tychonov cube $I^\kappa$. The image $f(B)$, being a compact subspace of the submetrizable space $X$, is metrizable. By \cite[2.7.12]{En}, the function $f$ depends on countably many coordinates, which means that there is a countable subset $C\subset \kappa$ such that $f=f_C\circ \pr_C$ where $\pr_C:I^\kappa\to I^C$ is the projection onto the  face $I^C$ of the cube $I^\kappa$ and $f_C:\pr_C(B)\to f(B)\subset X$ is a suitable continuous map. Since $X$ is an ANE for compact metrizable spaces, the map $f_C$ has a continuous extension $\tilde f_C:U\to X$ defined on an open neighborhood $U$ of $\pr_C(B)$ in the cube $I^C$. It follows that $V=\pr_C^{-1}(U)\cap A$ is an open neighborhood of $B$ in $A$ and $\tilde f=\tilde f_C\circ \pr_C|V:V\to X$ is a continuous extension of the map $f$, witnessing that $X$ is an ANE for the class of compact Hausdorff spaces.
\smallskip

Next, we show that $X$ is an ANE for the class of $k_\w$-spaces. Let $f:B\to X$ be a continuous map defined on a closed subset $B$ of a $k_\w$-space $A$. Then $A=\dlim A_n$ for some $k_\w$-sequence $(A_n)_{n\in\w}$ of compact subsets of $A$. Let $A_{-1}=\emptyset$. By induction, for each $n\in\w$ we can construct a continuous map $f_n:N_n(A_n\cap B)\to X$ defined on a closed neighborhood $N(B\cap A_n)$ of $B\cap A_n$ in $A_n$ and such that
\begin{itemize}
\item $N_n(B\cap A_n)\supset N_{n-1}(B\cap A_{n-1})$,
\item $f_n|B\cap A_n=f|B\cap A_n$ and 
\item $f_n|N_{n-1}(B\cap A_n)=f_{n-1}$. 
\end{itemize}
The inductive step can be done because $X$ is an ANE for the class of compact Hausdorff spaces. After completing the inductive construction, consider the set $N(B)=\bigcup_{n\in\w}N_n(B\cap A_n)$ and a map $\tilde f=\bigcup_{n\in\w}f_n:N(B)\to X$, which is a desired continuous extension of $f$ onto the open neighborhood $N(B)$ of $B$ in $A$.
\smallskip

The implication $(2)\Ra(1)$ trivially follows from the definitions of an $\ANRsw$ and $\mathsf{ANE}(k_\w)$-spaces.
\smallskip

Now assume that $X=\dlim X_n$ is the direct limit of a $k_\w$-sequence $(X_n)_{n\in\w}$ consisting of compact ANR's. We claim that $X$ is an ANE for the class of compact metrizable spaces. Let $f:B\to X$ be a continuous map defined on a closed subspace $B$ of a compact metrizable space $A$. Since $X$ carries the direct limit topology $\dlim X_n$, the compact subset $f(B)$ lies in some set $X_n$, $n\in\w$. Since $X_n$ is an ANR, the map $f:X\to X_n$ has a continuous extension $\tilde f:N(B)\to X_n\subset X$ defined on a neighborhood $N(B)$ of $B$ in $A$. 
\end{proof}

\section{Some subfunctors of the functor $F_\K$}

In the proof of Theorem~\ref{main} we shall apply a deep Basmanov's result on the preservation of compact ANR's by monomorphic functors of finite degree in the category $\Comp$ of compact Hausdorff spaces and their continuous maps. Let $\C$ be a full subcategory of the category $\Top$, containing all finite discrete spaces. 

We say that a functor $F:\C\to\Top$  
\begin{itemize}
\item is {\em monomorphic} if $F$ preserves monomorphisms (which coincide with injective continuous maps in the category $\Top$ and its full subcategory $\C$);
\item has {\em finite supports} ({\em degree} $\deg F\le n$) if for each object $X$ of the category $\C$ and each element $a\in FX$ there is a map $f:A\to X$ of a finite discrete space $A$ (of cardinality $|A|\le n$) such that $a\in Ff(FA)$;
\end{itemize}
The smallest number $n\in\w$ such that $\deg F\le n$ is called the {\em degree} of $F$ and is denoted by $\deg F$. If no such a number $n\in\w$ exists, then we put $\deg F=\infty$.

The following improvement of the classical Basmanov's theorem \cite{Bas} was recently  proved in \cite{BKZ}.

\begin{theorem}\label{BKZ} Let $F:\Comp\to\Comp$ be a monomorphic functor of finite degree $n=\deg F$ such that the space $Fn$ is finite. Then 
the functor $F$ preserves the classes of metrizable, finite-dimensional, and ANR-compacta.
\end{theorem} 

We shall apply this theorem to the subfunctors $\langle\cdot\rangle^L_n$ of the functor $F_\K$.
We recall that $\K$ is a non-trivial complete quasivariety of topological $E$-algebras of countable discrete signature $E$. By Theorem~\ref{choban2}, $F_\K$ can be thought as a functor $F_\K:\Kw\to\Kw$ in the category $\Kw$ of $k_\w$-spaces and their continuous maps.
By Theorem~2.4 of \cite{Chob2}, for each Tychonov space $X$ the free topological $E$-algebra $F_\K(X)$ is algebraically free in the sense that any bijective map $i:X_d\to X$ from a discrete topological space $X_d$ induces an algebraic isomorphism $F_\K i:F_\K X_d\to F_\K X$. This fact implies:

\begin{lemma}\label{l1} The functor $F_\K:\Tych\to\Top$ is monomorphic.
\end{lemma}

\begin{proof} Let $f:X\to Y$ be an injective continuous map between Tychonov  
spaces and $f_d:X_d\to Y_d$ be the same map between these spaces endowed with the discrete topologies. Let $i_X:X_d\to X$ and $i_Y:Y_d\to Y$ be the identity maps. 
Let be any (automatically continuous) map $r:Y_d\to X_d$ such that $r\circ f_d=\id_{X_d}$. Thus we obtain the commutative diagram:
$$\xymatrix{
X\ar[r]^f&Y\\
X_d\ar[u]^{i_X}\ar@<2pt>[r]^{f_d}&Y\ar@<2pt>[l]^r\ar[u]_{i_Y}
}
$$

Applying the functor $F_\K$ to this diagram we get the diagram
$$\xymatrix{
F_\K X\ar[r]^{F_\K f}&F_\K Y\\
F_\K X_d\ar[u]^{F_\K i_X}\ar@<2pt>[r]^{F_\K f_d}&F_\K Y\ar@<2pt>[l]^{F_\K r}\ar[u]_{F_\K i_Y}
}
$$The ``vertical'' maps $F_\K i_X:F_\K X_d\to F_\K X$ and $F_\K i_Y:F_\K Y_d\to F_\K Y$ in this diagram are bijective because the algebras $F_\K X$ and $F_\K Y$ are algebraically free. 
Taking into account that $F_\K r\circ F_\K f_d=F_\K(r\circ f_d)=F_\K\id_{X_d}=\id_{F_\K X_d}$, we conclude that the map $F_\K f_d$ is injective and so is the map $F_\K f:F_\K X\to F_\K Y$ becuase of the bjectivity of the maps $F_\K i_X$ and $F_\K i_Y$.
\end{proof}


 Now for every compact subset $L\subset E$ and every $n\in\w$ consider the functor $\langle \cdot\rangle_n^L:\Comp\to\Comp$ which assigns to each compact Hausdorff space $X$ the subspace $\langle X\rangle_n^L$ of $F_\K X$. The functor $\langle \cdot\rangle^L_n$ assigns to each continuous map $f:X\to Y$ between compact Hausdorff spaces the restriction $\langle f\rangle^L_n=F_\K f|\langle X\rangle^L_n$ of the homomorphism $F_\K f:F_\K X\to F_\K Y$.  

\begin{lemma}\label{l2} For every $n\in\IN$, $\langle\cdot\rangle^L_n:\Comp\to\Comp$ is a well-defined monomorphic functor of finite degree in the category $\Comp$.
\end{lemma}

\begin{proof} First we check that for each continuous map $f:X\to Y$ between compact Hausdorff spaces, the morphism $\langle f\rangle^L_n=F_\K f|\langle X\rangle^L_n$ is well-defined, which means that $F_\K f(\langle X\rangle^L_n)\subset \langle Y\rangle^L_n$. This will be done by induction on $n\in\w$. 

For $n=0$ the inclusion $F_\K(\langle X\rangle^L_0)=F_\K(X)=f(X)\subset Y=\langle Y\rangle^L_0$ follows from the fact that the homomorphism $F_\K$ extends the map $f$ (here we identify $X$ and $Y$ with the subspaces $\eta(X)$ and $\eta(Y)$ in $F_\K(X)$ and $F_\K(Y)$, respectively).

Assume that the inclusion $F_\K f(\langle X\rangle^L_{n})\subset \langle Y\rangle^L_{n}$ has been proved for some $n\in\w$. By definition,
$$\langle X\rangle^L_{n+1}=\langle X\rangle^L_{n}\cup\bigcup_{k\in\w}e_{k,X}((E_k\cap L)\times (\langle X\rangle^L_{n})^k).$$ Fix any element $x\in \langle X\rangle^L_{n+1}$. If $x\in \langle X\rangle^L_{n}$, then $$F_\K(x)\in F_\K(\langle X\rangle^L_n)\subset \langle Y\rangle^L_{n}\subset\langle Y\rangle^L_{n+1}$$ by the inductive assumption.

If  $x\in \langle X\rangle^L_{n+1}\setminus \langle X\rangle^L_{n}$, then $x=e_{k,X}(z,x_1,\dots,x_k)$ for some $k\in\w$, $z\in E_k\cap L$, and points 
$x_1,\dots,x_k\in \langle X\rangle^L_{n}$. Since $F_\K f$ is a $E$-homomorphism, we get
$$
\begin{aligned}
F_\K f(x)&=F_\K f(e_{k,X}(z,x_1,\dots,x_k))=e_{k,Y}(z,F_\K f(x_1),\dots,F_\K f(x_k))\in\\
& \in e_{k,Y}((E_k\cap L)\times (\langle Y\rangle^L_{n})^k)\subset\langle Y\rangle^L_{n+1}.
\end{aligned}$$
Thus for every $n\in\w$ the functor $\langle \cdot\rangle^L_n$ is well-defined. It is monomorphic as a subfunctor of the monomorphic functor $F_\K$.

Next, we show that the functor $\langle \cdot \rangle$ has finite degree. This will be done by induction on $n\in\w$. Since $\langle X\rangle^L_{0}=X$, $\deg \langle\cdot\rangle^L_{0}=1$.

Assume that for some $n\in\w$ the functor $\langle\cdot\rangle^L_{n}$ has finite deree $d$. Since $L$ is a compact subset of $E$, there is $m\in\w$ such that $L\cap E_k=\emptyset$ for all $k\ge m$. We claim that $\deg\langle\cdot\rangle^L_{n+1}\le m\cdot d$. Take any element $x\in \langle X\rangle^L_{n+1}$. If $x\in\langle X\rangle^L_{n}$, then by the inductive assumption there is a subset $A\subset X$ of cardinality $|A|\le d$ such that $x\in\langle A\rangle^L_{n}$ and we are done. If $x\in \langle X\rangle^L_{n+1}\setminus\langle X\rangle^L_{n}$, then $x=e_{k,X}(z,x_1,\dots,x_k)$ for some $k\in\w$, $z\in E_k\cap L$, and points $x_1,\dots,x_k\in\langle X\rangle^L_{n}$. Since $L\cap E_k\ni z$ is not empty, $k\le m$. By the inductive assumption, for every $i\le k$ there is a finite subset $A_i\subset X$ of cardinality $|A_i|\le d$ such that $x_i\in\langle A_i\rangle^L_{n}$. Then the union $A=\bigcup_{i=1}^kA_i$ has cardinality $|A|\le k\cdot d\le m\cdot d$ and 
$$x=e_{k,X}(z,x_1,\dots,x_k)\in e_{k,X}((L\cap E_k)\times(\langle A\rangle^L_{n})^k)\subset\langle A\rangle^L_{n+1}$$witnessing that the functor $\langle\cdot\rangle^L_{n+1}$ has finite degree $\deg \langle \cdot\rangle^L_{n+1}\le m\cdot d$. 
\end{proof}

\begin{lemma}\label{l3} If $L\subset E$ is finite, then for each $n\in\w$ the functor $\langle\cdot \rangle^L_n$ preseves finite spaces.
\end{lemma}

\begin{proof} Let $X$ be a finite space. By induction on $n\in\w$ we shall show that the space $\langle X\rangle^L_n$ is finite. This is clear for $n=0$. Assume that for some $n\in\w$ the space $\langle X\rangle^L_{n}$ is finite. Since $L\subset E$ is finite there is $m\in\w$ such that $L\cap E_n=\emptyset$ for all $k>m$. Then $$\langle X\rangle^L_{n+1}=\langle X\rangle^L_{n}\cup\bigcup_{k\le m}e_{k,X}((E_k\cap L)\times(\langle X\rangle^L_{n})^k)$$is finite as the finite union of finite sets.
\end{proof}

Combining Lemmas~\ref{l2}, \ref{l3} with Theorem~\ref{BKZ}, we get

\begin{corollary}\label{corANR} For any finite subset $L\subset E$ and every $n\in\w$ the functor $\langle\cdot\rangle^L_n$ preserves (finite-dimensional) compact ANR's.
\end{corollary}

\section{Proof of Theorem~\ref{main}}

Without loss of generality, the quasivariety $\K$ is non-trivial (otherwise, $F_\K(X)$ is a singleton and hence is an $\ANRsw$-space for each non-empty space $X$).

Let $X$ be a submetrizable $\ANRsw$-space. By the ANR-Theorem \ref{ANR}, the product $X\times Q^\infty$ is a $Q^\infty$-manifold. By the Triangulation Theorem \ref{triangle}, $X\times Q^\infty$ is homeomorphic to $T\times Q^\infty$ for a  countable locally finite simplicial complex $T$. This implies that $X\times Q^\infty$ can be written as the direct limit $X\times Q^\infty=\dlim X_n$ of a $k_\w$-sequence $(X_n)_{n\in\w}$ of compact ANR's. 

Write the countable discrete space $E$ as the direct limit $E=\dlim L_n$ of a $k_\w$-sequence $(L_n)_{n\in\w}$ 
of finite subsets of $E$. By Choban's Theorem~\ref{choban2}, the space $F_\K(X\times Q^\infty)$ is the direct limit $\dlim \langle X_n\rangle_n^{L_n}$ of the $k_\w$-sequence $\langle X_n\rangle^{L_n}_n$. By Corollary~\ref{corANR}, each space $\langle X_n\rangle^{L_n}_n$, $n\in\w$, is a compact metrizable ANR. Consequently, $F_\K (X\times Q^\infty)=\dlim \langle X_n\rangle^{L_n}_n$ is a submetrizable $\ANRsw$-space by Theorem~\ref{ANRkw}. Since $X$ is retract of $X\times Q^\infty$, the space $F_\K X$ is a retract of $F_\K(X\times Q^\infty)$ and hence $F_\K X$ is a submetrizable $\ANRsw$-space.

Now assume that $X$ is a  compactly finite-dimensional $s_\w$-space. 
Then $X=\dlim X_n$ is the direct limit of finite-dimensional compact metrizable spaces. By the Choban's Theorem~\ref{choban2}, the space $F_\K(X\times Q^\infty)$ is the direct limit $\dlim \langle X_n\rangle_n^{L_n}$ of the $k_\w$-sequence $\langle X_n\rangle^{L_n}_n$. Corollary~\ref{corANR} implies that each compact space $\langle X_n\rangle^{L_n}_n$ is  metrizable and finite-dimensional. Then the space $F_\K X=\dlim \langle X_n\rangle^{L_n}_n$ is compactly finite-dimensional, being  the direct limts of finite-dimensional compacta.

\end{document}